\documentclass[12pt]{article}
\usepackage{amsxtra,amssymb,amsthm,amsmath,latexsym}

\theoremstyle{plain}
\newcommand{\refS}[1]{Section~\ref{S:#1}}

\def\oH{{\overset{\circ}{H}}}
\def\oH1{{\overset{\circ}{H}\kern-.02in{}^1}}

\def\bee{\begin{equation*}}
\def\eee{\end{equation*}}
\def\be{\begin{equation}}
\def\ee{\end{equation}}

\begin{document}

%\begin{titlepage}
\title{ Existence and uniqueness of the global solution to the Navier-Stokes equations
}

\author{Alexander G. Ramm\\
 Department  of Mathematics, Kansas State University, \\
 Manhattan, KS 66506, USA\\
ramm@math.ksu.edu\\
%,\\ %fax 785-532-0546, tel. 785-532-0580}
%http://www.math.ksu.edu/\,$\widetilde{\ }$\,ramm}
}

\date{}
\maketitle\thispagestyle{empty}

%%%%%%%%%%%%%%%%%%%%%%%%%%%%%%%%%%%%%%%%%%%%%%%%%
\begin{abstract}

A proof is given of the global existence and uniqueness of a weak  solution to Navier-Stokes
equations in unbounded exterior domains.
\end{abstract}
%%%%%%%%%%%%%%%%%%%%%%%%%%%%%%%%%%%%%%%%%%%%%%%%%
%\end{titlepage}

\section{Introduction}\label{S:1}
Let $D\subset \mathbb{R}^3$ be a bounded domain with a connected $C^2-$smooth boundary $S$,
and $D':=\mathbb{R}^3\setminus D$ be the unbounded exterior domain.

Consider the Navier-Stokes equations:
\be\label{e1} u_t+(u, \nabla)u=-\nabla p +\nu\Delta u +f, \quad x\in D',\,\, t\ge 0,
\ee
\be\label{e2} \nabla \cdot u=0,
\ee
\be\label{e3} u|_S=0, \quad u|_{t=0}=u_0(x).
\ee
Here $f$ is a given vector-function, $p$ is the pressure, $u=u(x,t)$ is the velocity vector-function,
$\nu=const >0$ is the viscosity coefficient, $u_0$ is the given initial velocity, $u_t:=\partial_t u$,
$(u, \nabla)u:= u_a \partial_a u$, $\partial_a u:=\frac {\partial u}{\partial x_a}:=u_{;a}$,
 and $ \nabla \cdot u_0:=u_{a;a}=0$.  Over the repeated indices $a$ and $b$ summation is
understood, $1\le a,b \le 3$.  All functions are assumed real-valued.

  We assume that $u\in W$,
$$W:=\{u| L^2(0,T; H^1_0(D'))\cap
L^\infty(0,T; L^2(D'))\cap u_t\in L^2(D'\times [0,T]); \nabla \cdot u=0 \},$$
 where  $T>0$ is arbitrary.

Let $(u,v):=\int_{D'} u_av_adx$ denote the inner product in $L^2(D')$,  $\|u\|:=(u,u)^{1/2}$. By $u_{ja}$ the $a-$th component of the vector-function $u_j$
is denoted, and $u_{ja;b}$ is the derivative $\frac{\partial u_{ja}}{\partial x_b}$.  Equation (2) can be written as $u_{a;a}=0$ in these notations. We denote $\frac{\partial u^2}{\partial x_a}:= (u^2)_{;a}$, $u^2:=u_bu_b$. By $c>0$ various estimation constants are
denoted.

{\em Let us define a weak solution to
problem (1)-(3) as an element of $W$ which satisfies the identity:}
\be\label{e4}
(u_t, v)+(u_au_{b;a}, v_b)+\nu(\nabla u, \nabla v)=(f,v),  \qquad \forall v\in W.
\ee

 Here we took into account that $-(\Delta u,v)= (\nabla u, \nabla v)$ and $(\nabla p, v)=-(p,v_{a;a})=0$
 if $v\in H^1_0(D')$ and $\nabla \cdot v=0$.
 Equation (4) is equivalent to the integrated equation:
$$
\int_0^t[(u_s, v)+(u_{a}u_{b;a}, v_{b})+\nu(\nabla u, \nabla v)]ds=\int_0^t(f,v)ds,  \qquad \forall v\in W \quad (*).
$$
Equation (4) implies equation (*), and differentiating equation (*) with respect to $t$ one gets equation (4) for almost all $t\ge 0$.

The aim of this paper is to prove the global existence and uniqueness  of the
weak solution to the Navier-Stokes boundary problem, that is, solution in $W$  existing for all $t\ge 0$.
Let us assume that
$$\sup_{t\ge 0}\int_0^t\|f\|ds\le c, \qquad (u_0,u_0)\le c. \qquad (A)$$

{\bf Theorem 1.}  {\em If assumptions (A) hold and $u_0\in H^1_0(D)$ satisfies equation (2),
then there exists for all $t>0$ a solution $u\in W$ to (4)  and this solution is unique in $W$
provided that $\|\nabla u\|^4\in L^1_{loc}(0, \infty)$.}
\vspace{3mm}

In \refS{2} we prove Theorem 1. There is a large literature on Navier-Stokes equations, of which we
mention only \cite{L} and \cite{T}. The global existence and uniqueness of the solution to
Navier-Stokes boundary problems has not yet been proved without additional assumptions.
Our additional assumption is  $\|\nabla u\|^4\in L^1_{loc}(0, \infty)$.
The history of this problem see, for example, in \cite{L}.
In \cite{T} the uniqueness of the global solution to Navier-Stokes equations is established under the assumption
$\|u\|^8_{L^4(D')}\in L^1_{loc}(0,\infty)$.

\section{ Proof of Theorem 1}\label{S:2}

\begin{proof}[Proof of Theorem 1] The steps of the proof are: a) derivation of a priori estimates; b) proof of the existence of the
solution in $W$; c) proof of the uniqueness of the solution in $W$.

a) {\em Derivation of a priori estimates.}

Take $v=u$ in (4). Then
$$(u_au_{b;a}, u_b)=-(u_a u_b, u_{b;a})=-\frac 1 2(u_a, (u^2)_{;a})= \frac 1 2(u_{a;a},u^2)=0,$$
where the equation $u_{a;a}=0$ was used. Thus, equation (4) with $v=u$ implies
\be\label{e5} \frac 1 2 \partial_t (u,u)+\nu (\nabla u, \nabla u)=(f,u)\le \|f\| \|u\|.
\ee
 We will use the known inequality $||u||||f||\le \epsilon ||u||^2+ \frac 1{4\epsilon}||f||^2$ with a small $\epsilon>0$, and denote by $c>0$  {\em various} estimation constants.

  One gets from (5) the following estimate:
 \be\label{e6}   (u(t),u(t))+ 2 \nu\int_0^t(\nabla u, \nabla u)ds\le (u_0,u_0)+2\int_0^t\|f\|ds \sup_{s\in [0,t]}\|u(s)\|\le c+c \sup_{s\in [0,t]}\|u(s)\|.
\ee
Recall that assumptions (A) hold.
Denote $\sup_{s\in [0,t]}\|u(s)\|:=b(t)$. Then inequality (6) implies
\be\label{e7}
b^2(t)\le c+cb(t), \qquad c=const >0.
\ee
Since $b(t)\ge 0$, inequality (7) implies
\be\label{e8}
\sup_{t\ge 0}b(t)\le c.
\ee
Remember that $c>0$ denotes various constants, and the constant in equation (8) differs from the constant
in equation (7).
From (8) and (6) one obtains
\be\label{e9}
\sup_{t\ge 0}[(u(t),u(t)) +  \nu\int_0^t(\nabla u, \nabla u)ds]\le c.
\ee
A priori estimate (9) implies for every $T\in [0,\infty)$ the inclusions
$$ u\in L^\infty(0,T; L^2(D')), \qquad u\in  L^2(0,T; H^1_0(D')).$$
This and equation (4) imply that $u_t\in L^2(D'\times [0,T])$ because equation (4)
shows that $(u_t,v)$ is bounded for every $v\in W$. Note that
$L^\infty(0,T; L^2(D'))\subset  L^2(0,T; L^2(D'))$, and that bounded sets
in a Hilbert space are weakly compact. Weak convergence is denoted by the sign
$\rightharpoonup$.

b) {\em Proof of the existence of the solution $u\in W$ to (4) and (*).}

The idea of the proof is to reduce the problem to the existence of the solution to
a Cauchy problem for ordinary differential equations (ODE) of finite order, and then to use a priori estimates to establish convergence of these solutions of ODE
to a solution of equations (4) and (*). This idea is used, for example, in  \cite{L}.
Our argument differs from the arguments in the literature in treating the limit of the term
$\int_0^t(u^n_s,v)ds$.

Let us look for a solution to equation (4) of the form $u^n:=\sum_{j=1}^n c_j^n(t)\phi_j(x)$,
where $\{\phi_j\}_{j=1}^\infty$ is an orthonormal basis of the space $L^2(D')$ of divergence-free
vector functions belonging to $H^1_0(D')$ and in the expression  $u^n$ the upper index $n$
is not a power. If one substitutes $u^n$ into equation (4),
takes $v=\phi_m$, and uses the orthonormality of the system $\{\phi_j\}_{j=1}^\infty$ and the relation
$(\nabla \phi_j, \nabla \phi_m)=\lambda_m \delta_{jm}$, where $\lambda_m$ are the eigenvalues of the
vector Dirichlet Laplacian in $D$ on the divergence-free vector fields, then one gets a system of ODE
for the unknown coefficients $c_m^n$:
\be\label{e10}
\partial_t c_m^n +\nu\lambda_m c_m^n +\sum_{i,j=1}^n (\phi_{ia}\phi_{jb;a}, \phi_{mb}) c_i^nc_j^n=f_m,\quad c_m^n(0)=(u_0,\phi_m).
\ee

Problem (10) has a unique global solution  because of the a priori estimate that follows from (9) and from Parseval's relations:
\be\label{e11}
\sup_{t\ge 0} (u^n(t),u^n(t))=\sup_{t\ge 0}\sum_{j=1}^n [c_j^n(t)]^2\le c.
\ee
  Consider the set $\{u^n=u^n(t)\}_{n=1}^\infty$.
Inequalities (9) and (11) for $u=u^n$ imply the existence
of the weak limits $u^n\rightharpoonup u$ in $ L^2(0,T; H^1_0(D'))$ and in $ L^\infty(0,T; L^2(D'))$.
This allows one to pass to the limit in equation (*) in all the terms except the first,
 namely, in the term  $\int_0^t(u^n_{s},v(s))ds$. The weak limit of the term $(u^n_{a}u^n_{b;a}, v_b)$
 exists and is equal to
 $(u_au_{b;a}, v_b)$ because
 $$(u^n_au^n_{b;a}, v_b)=-(u^n_{a}u^n_{b}, v_{b;a})\to -(u_au_b,v_{b;a})=(u_au_{b;a},v_b).$$
 Note that  $ v_{b;a}\in L^2(D')$ and  $u^n_{a} u^n_{b}\in L^4(D')$.
The relation $(u^n_au^n_{b;a}, v_b)=-(u^n_{a}u^n_{b}, v_{b;a})$
follows from an integration by parts and from the equation $u^n_{a;a}=0$.

 The following inequality is essentially known:
\be\label{e12}
\|u\|_{L^4(D')}\le 2^{1/2}\|u\|^{1/4}\|\nabla u\|^{3/4}, \quad \|u\|:=\|u\|_{L^2(D')}, \quad u\in H^1_0(D').
\ee
In \cite{L} this inequality is proved for $D'=\mathbb{R}^3$, but a function $u\in  H^1_0(D')$
can be extended by zero to $D=\mathbb{R}^3\setminus D'$ and becomes an element of $H^1( \mathbb{R}^3)$ to which
inequality (12) is applicable.

It follows from (12) and the Young's inequality ($ab\le \frac {a^p}{p}+\frac {b^q}{q}$, $p^{-1}+q^{-1}=1$) that
\be\label{e13}
\|u\|^2_{L^4(D')}\le \epsilon \|\nabla u\|^{2}+\frac{27}{16 \epsilon^3}  \|u\|^2, \quad u\in H^1_0(D'),
\ee
where $\epsilon>0$ is an arbitrary small number, $p=\frac 4 3$ and $q=4$.  One has
$u^n_a u^n_{b}\rightharpoonup u_au_b$ in $L^2(D')$  as $n\to \infty$, because bounded sets
in a reflexive Banach space $L^4(D')$ are weakly compact.
Consequently, $(u^n_au^n_{b;a}, v_b)\to (u_au_{b;a}, v_b)$ when $n\to \infty$, as claimed.
Therefore, $ \int_0^t(u^n_au^n_{b;a}, v_b)ds\to \int_0^t(u_au_{b;a}, v_b)ds$. The weak limit
of the term $\nu \int_0^t(\nabla u^n, \nabla v)ds$ exists because of the a priori estimate (9)
and the weak compactness of the bounded sets in a Hilbert space.
Since equation (*) holds, and the limits of all its terms, except $\int_0^t(u^n_s,v)ds$,
do exist, then  there exists the limit $\int_0^t (u^n_{s},v(s))ds\to \int_0^t(u_s, v(s))ds$
 for all $v\in W$. By passing to the limit $n\to \infty$
 one proves that the limit $u$ satisfies equation (*). Differentiating equation (*) with respect to $t$ yields equation (4) almost everywhere.

c) {\em Proof of the uniqueness of the solution $u\in W$}.

 Suppose there are two solutions to equation (4), $u$ and $w$,  $u,w\in W$,
and let $z:=u-w$. Then
\be\label{e14}
(z_t, v) +\nu (\nabla z, \nabla v) +(u_au_{b;a}- w_aw_{b;a}, v_b)=0.
\ee
Since $z\in W$, one may set $v=z$ in (14) and get
\be\label{e15}
(z_t, z) +\nu (\nabla z, \nabla z) + (u_au_{b;a}- w_aw_{b;a}, z_b)=0, \qquad z=u-w.
\ee
Note that $(u_au_{b;a}- w_aw_{b;a}, z_b)=(z_au_{b;a},z_b)+(w_a z_{b;a},z_b)$,
and $(w_a z_{b;a},z_b)=0$ due to the equation $w_{a;a}=0$. Thus, equation (15) implies
\be\label{e16}
\partial_t (z, z) +2\nu (\nabla z, \nabla z)\le 2|(z_au_{b;a}, z_b)|.
\ee
Since $|z_a  u_{b;a} z_b|\le |z|^2 |\nabla u|$, one has the following estimate:
\be\label{e17}
|(z_au_{b;a}, z_b)|\le \int_{D'} |z|^2 |\nabla u|dx\le \|z\|^2_{L^4(D')} \|\nabla u\|\le
 \|\nabla u\|\Big(\epsilon \|\nabla z\|^2+ \frac{27}{16\epsilon^3}\|z\|^2\Big).
\ee
 Denote $\phi:=(z,z)$, take into account that $\|\nabla u\|^4\in L^1_{loc}(0,\infty)$,
 choose $\epsilon=\frac {\nu}{ \|\nabla u\|}$ in the inequality (13), in which $u$ is replaced by $z$,
  use inequality (17) and get
\be\label{e18}
\partial_t \phi+\nu (\nabla z, \nabla z) \le \frac{27}{16\nu^3}\|\nabla u\|^4 \phi, \qquad \phi|_{t=0}=0.
\ee
In the derivation of inequality (18) the idea is to compensate the term $\nu \|\nabla z\|^2$ on the left
side of inequality (16)
by the term $\epsilon \|\nabla u\| \|\nabla z\|^2$ on the right side of inequality (17).
To do this, choose $\|\nabla u\|\epsilon=\nu$ and obtain inequality (18).
It follows from inequality (18) that
$$\partial_t\phi \le \frac {27 \|\nabla u\|^4}{16\nu^3}  \phi, \quad \phi|_{t=0}=0.$$
Since we have assumed that $\|\nabla u\|^4\in L^1_{loc}(0,\infty)$ this implies that $\phi=0$ for all $t\ge 0$.

Theorem 1 is proved.
\end{proof}

{\bf Remark 1.} One has (summation is understood over the repeated indices):
$$2|(z_au_{b;a},z_b)|=2|(z_au_b, z_{b;a})|\le 18\|\nabla z\|\||z||u|\|\le \nu \|\nabla z\|^2 +\frac {81} {\nu}\||z||u|\|^2.$$
Thus,
$$ \partial_t \phi+\nu (\nabla z, \nabla z) \le \frac {81} {\nu}\||z||u|\|^2.$$
If one assumes that $|u(\cdot,t)|\le c(T) $ for every $t\in [0,T]$, then
$ \partial_t \phi\le c\phi, \quad \phi(0)=0,$ on any interval $[0,T]$,
 $c=c(T,\nu)>0$ is a constant. This implies $\phi=0$ for all $t\ge 0$. The same conclusion
 holds under a weaker assumption $\|u(\cdot,t)\|_{L^4(D')}\le c(T)$ for  every $t\in [0,T]$,
  or under even weaker assumption $\|u(\cdot,t)\|^8_{L^4(D')}\in L^1_{loc}(0,\infty)$.
\vspace{2mm}

In \cite{L} it is shown that the smoothness properties of the solution $u$
are improved when the smoothness properties of $f$, $u_0$ and $S$ are improved.

%\newpage


\begin{thebibliography}{1000} %number of characters of longest bibitem label

\bibitem{L} O. Ladyzhenskaya, The mathematical theory of viscous incompressible flow,
Gordon and Breach, New York, 1969.
%\bibitem{S} H.Sohr, The  Navier-Stokes equations, Birkh\"auser, Basel, 2001.

%\bibitem{R} A.G.Ramm, Large-time behavior of the weak solution to 3D Navier-Stokes equations,
%Appl. Math. Lett., 26 (2013), 252-257.

\bibitem{T} R.Temam, Navier-Stokes equations. Theory and numerical analysis,
North Holland, Amsterdam, 1984.



\end{thebibliography}
\end{document}